\documentclass{article}
\usepackage{amssymb}
\usepackage{graphicx}
\usepackage[latin1]{inputenc}
\newtheorem{theorem}{Theorem}
\def\biglf{\par\bigskip\noindent}
\newcommand{\bql}[1]{%
\begin{equation}\label{#1}%
}
\newcommand{\eq}{\end{equation}}
\newcommand{\R}{\ensuremath{\mathbb{R}}}
\def\qed{\hfill$\Box$}
\def\eref#1{(\ref{#1})}%
\newcommand{\MS}{\hbox{\tt sum}}
\begin{document}
\begin{center}
  {\bf
    A Meshfree
    Method for Solving the Monge--Amp\`ere Equation}

  Klaus B\"ohmer\footnote{Fachbereich Mathematik und Informatik,
Universit\"at Marburg,
Arbeitsgruppe Numerik,
Hans Meerwein Stra\ss{}e, Lahnberge,
D-35032 Marburg,
Germany},
  Robert Schaback \footnote{Institut f\"ur Numerische und Angewandte Mathematik,
Universit\"at G\"ottingen,
Lotzestra\ss{}e 16-18,
D-37073 G\"ottingen,
Germany}\\

~

\end{center}

{\bf Abstract}: 
This paper solves the two-dimensional Dirichlet problem for the
Monge-Amp\`ere equation by a strong meshless collocation
technique that uses a polynomial trial space and
collocation in the domain and
on the boundary. Convergence rates may be up to exponential,
depending on the smoothness of the true solution,
and this is demonstrated numerically and proven theoretically,
applying a sufficiently fine collocation discretization.
A much more thorough investigation
of meshless methods for fully nonlinear problems is in preparation.\\~\\

~\\


%%%%%%%%%%%%%%%%%%%%%
%****************************************************************
\section{Introduction}\label{SecIntro}
The Dirichlet problem for the Monge-Amp\`ere equation
on a bounded domain
${\Omega}\subset\R^d$
consists in finding a smooth function $u$ on 
$\overline{\Omega}$
such that the equations 
\bql{eqMA}
\begin{array}{rcl}
u_{xx}u_{yy}-u_{xy}^2 &=& g \hbox{ in } \overline{\Omega}\\
u &=& f  \hbox{ in }\partial \Omega
\end{array}
\eq
hold, where the functions $f$ and $g$ are given on the boundary
$\partial\Omega$ and on 
$\overline{\Omega}$, respectively. The goal of this contribution is to show
how %\textbf{Klaus simply }
simply meshless methods in strong form can be applied
% \textbf{Klaus
and their convergence be proven for this most important  special case
of a fully nonlinear second-order
partial differential equation, thoroughly investigated in B\"ohmer/Schaback
\cite{boehmer-schaback:2016-1}. Concerning simplicity,
%\textbf{Klaus Die neuen Zitate sind nach end{document} im bibtex mode vangehaengt!!!
the meshless method proposed here % This
is strongly superior %completely different for 
to the complicated approaches for difference methods in
Oberman \cite{oberman:2008-1} and finite elements.
 The first method including %a complicated
   convergence %proof is % has been 
 is published in B{\"{o}}hmer \cite{boehmer:2004-1}, %. Essential  were the results for
 based upon  smooth finite elements of Davydov and Saeed
 \cite{davydov:2006-1,davydov-saaed:2013-1}. See also  
 Brenner et al. \cite{brenner-et-al:2011-1}
 and Feng/Neilan \cite{feng-neilan:2009-1}
 for a method to solve the Monge-Amp\`ere
equations
via finite elements and a $C^0$ penalty method or a vanishing moment method. 
% with completely different ideas starting  2009.  

%Davydov06DaSa13  Neilan??BrNe13
%I just noticed that I still lack some of  Davydov with Saeed and  Brenner with Neilan work.
\biglf
To explain the connection to the general situation,
the equation is rewritten as
$$
\begin{array}{rcll}
Fu&=& g& \hbox{ in }\overline{\Omega}\\
u&=&f & \hbox{ in }\Gamma:=\partial \Omega
\end{array} 
$$
with the nonlinear map $F\;:\;u\mapsto u_{xx}u_{yy}-u_{xy}^2$.  
which is  defined as a mapping 
\bql{eqLuf0}
F\;:\; C^2(\overline{\Omega}) \to C(\overline{\Omega}).
\eq
Fully nonlinear problems are those where second-order
derivatives arise nonlinearly in $F$.
Existence and regularity results from the literature
% books Gilbarg/Trudinger \cite{gilbarg-trudinger:2001-1},
% Skrypnik \cite{skrypnik:1986-1}, and Chen/Wu \cite{chen-wu:1998}
are collected in \cite{boehmer:2010-1}, Theorems 2.79-2.82.
\biglf
Under appropriate assumptions, there is a unique solution $u^*\in
C^2(\overline{\Omega})$. Then 
the connection to linear elliptic problems is made by the
linearization at $u^*$, namely
$$
F'(u^*)(v)=u^*_{xx}v_{yy}+v_{xx}u^*_{yy}-2u^*_{xy}v_{xy}.
$$
By the cited background literature collected in
% Gilbarg/ Trudinger \cite{gilbarg-trudinger:2001-1},
% Trudinger and Wang \cite{trudinger-wang:2008-1}  and B{\"o}hmer
\cite{boehmer:2010-1}, this $F'(u^* )$ defines an elliptic operator
for a convex  locally unique  solution $u^*$. 
There also is a variation of a maximum principle, and 
$F'(u^*)$
is boundedly invertible. %%%%%%%%%%%%%%%%%%%%%%%
%%%%%%%%%%%%%%%%%%%%%%%%%%%%%%%%%%%%%%%%%%%%%%%%%%
%****************************************************************
\section{Strong Meshfree Discretizations}\label{SecMD}
These techniques discretize PDE problems
using 
\begin{enumerate}
\item a {\em trial space} that
  should approximate the true solution $u^*$ well, and 
\item sets of {\em test points} on which the differential operator and the
  boundary conditions are directly sampled,
  \item forming a nonlinear system of collocation equations that is possibly
    overdetermined, and finally
    \item applying a nonlinear optimizer to minimize residuals of the system.
\end{enumerate}
In contrast to finite element methods,
there is no connection between the test and the trial side via a
triangularization, and there is no numerical integration.
{\em Consistency} is guaranteed by choosing a sufficiently rich
trial space, and {\em stability} requires to choose sufficiently many
well-posed collocation points. The details concerning these choices
are nontrivial and will not be explained here, see e.g.
\cite{schaback:2010-2,boehmer-schaback:2013-1,schaback:2015-4}
for a comprehensive convergence analysis.
\biglf
To present a simple example that works for the Monge-Amp\`ere equation,
we confine the domain $\Omega$ to be the unit square $[-1,+1]^2$  
and use polynomial trial functions
\bql{equSScxy}
u^C(x,y):=\displaystyle{\sum_{m=0}^M\sum_{n=0}^{M-m}c_{mn}x^my^n} 
\eq
of total degree at most $M$
that will approximate smooth solutions well. The coefficients
are collected into a triangular matrix $C$. Other parametrized
trial spaces could serve the same purpose, but the trial functions
should be smooth and derivatives up to second order should be
available at low computational cost.
\biglf
On the test side, we take
$K_D$
points $z_i^D:=(x_i^D,y_i^D),\;1\leq
i\leq K_D$ in the closure of the domain for
approximation of the
differential operator, and   additionally
$K_B$ points $z_i^B:=(x_i^B,y_i^B),\;1\leq
i\leq K_B$ 
on the boundary, forming point sets $Z^D$ and $Z^B$.
These sets may have a nonempty intersection on the boundary.
\biglf
The strong meshless discretization then sets up a nonlinear system 
$$
\begin{array}{rcl}
  u^C_{xx}(z_i^D)u^C_{yy}(z_i^D)-u^C_{xy}(z_i^D)^2 &=& g(z_i^D),\;
  1\leq i\leq K_D\\
  u^C(z_j^B) &=& f(z_j^B),\;    1\leq i\leq K_B
\end{array} 
$$
for the coefficient matrix $C$. There is
no linearization done here. Linearization is left to the solver.
\biglf
Brute-force numerical methods
can then apply either a nonlinear equation solver
to the system or a nonlinear optimizer to the residuals.
Because oversampling will often be necessary to guarantee stability
\cite{schaback:2010-2,boehmer-schaback:2013-1,schaback:2015-4},
the latter situation is preferable.
This paper uses the nonlinear least-squares minimizer
{\tt nlsqnonlin} of MATLAB on the
residuals of the above nonlinear system.
No matter how the nonlinear solver works, the possibility
of multiple local solutions requires good initial startup parameters. 
If none are known from external arguments, a standard
technique is to apply repeated calculations using larger and larger trial
spaces, starting from the optimal solution of the previous step.
In parallel, the test discretizations should be refined from step to step.
See Fasshauer's book
\cite{fasshauer-mccourt:2015-1} for a comprehensive account
of meshless methods using MATLAB.
\biglf
A similar meshless algorithm, but without a convergence analysis,
and with different trial functions,
was proposed by Zhiyong Liu and collaborators
\cite{liu-he:2013-1,liu-he:2014-1,li-liu:2017-1},
based on Kansa's unsymmetric collocation.
Finite differences and multigrid methods were applied by various authors
\cite{benamou-et-al:2010-1,froese-oberman:2011-1,
  froese-oberman:2013-1,liu-et-al:2017-1}, but these techniques
are further away
from this work because they are not meshless and much more difficult to
implement. 
\biglf
Readers may jump to Section \ref{secMATDet}
for MATLAB implementation details, and to Section \ref{SecConv}
for the theoretical error analysis. We shall present the
numerical results next, and finally add a summary and an outlook.
%****************************************************************
%****************************************************************
\section{Numerical Results}\label{SecR1}
\begin{figure}[hbt]
\begin{center}
\includegraphics[width=10cm]{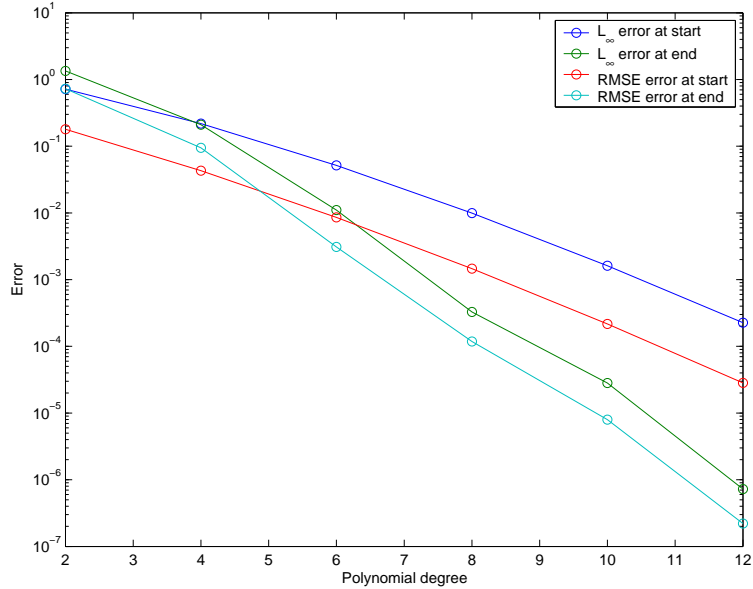}\\
\caption{Errors as functions of degree}\label{figtotal}
\end{center}
\end{figure}

\begin{figure}[hbt]
\begin{center}
\includegraphics[width=10cm]{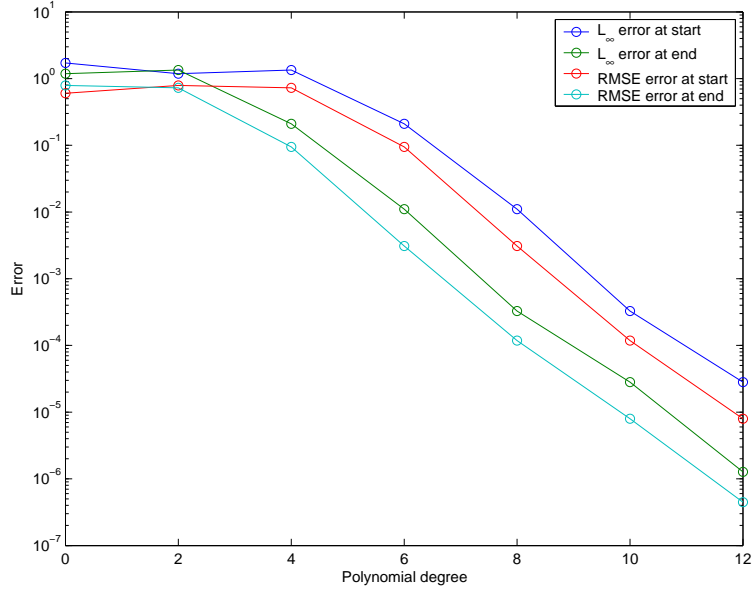}\\
\caption{Errors as functions of degree}\label{figtotal2}
\end{center}
\end{figure}

The example of this paper uses the true solution
\bql{equexp} 
u^*(x,y)=\exp((x^2+y^2)/2)
\eq
to generate the appropriate functions $f$ and $g$ in \eref{eqMA}.
\biglf
If the truncated Taylor expansion of the true solution is
used as a starting approximation, Figure
\ref{figtotal}
shows the exponential decay of the error as a function of the total degree
$M$ in \eref{equSScxy}. 
Here, RMSE stands for the root-mean-square error. 
If we start with the constant approximation 1 at degree zero and use the best
coefficients for degree $M$ to start at degree $M+2$, we get Figure \ref{figtotal2}.
We add figures for the final $M=12$ calculation, starting from the 
optimal $M=10$ result.
\begin{figure}[hbt]
\begin{center}
\includegraphics[width=10cm]{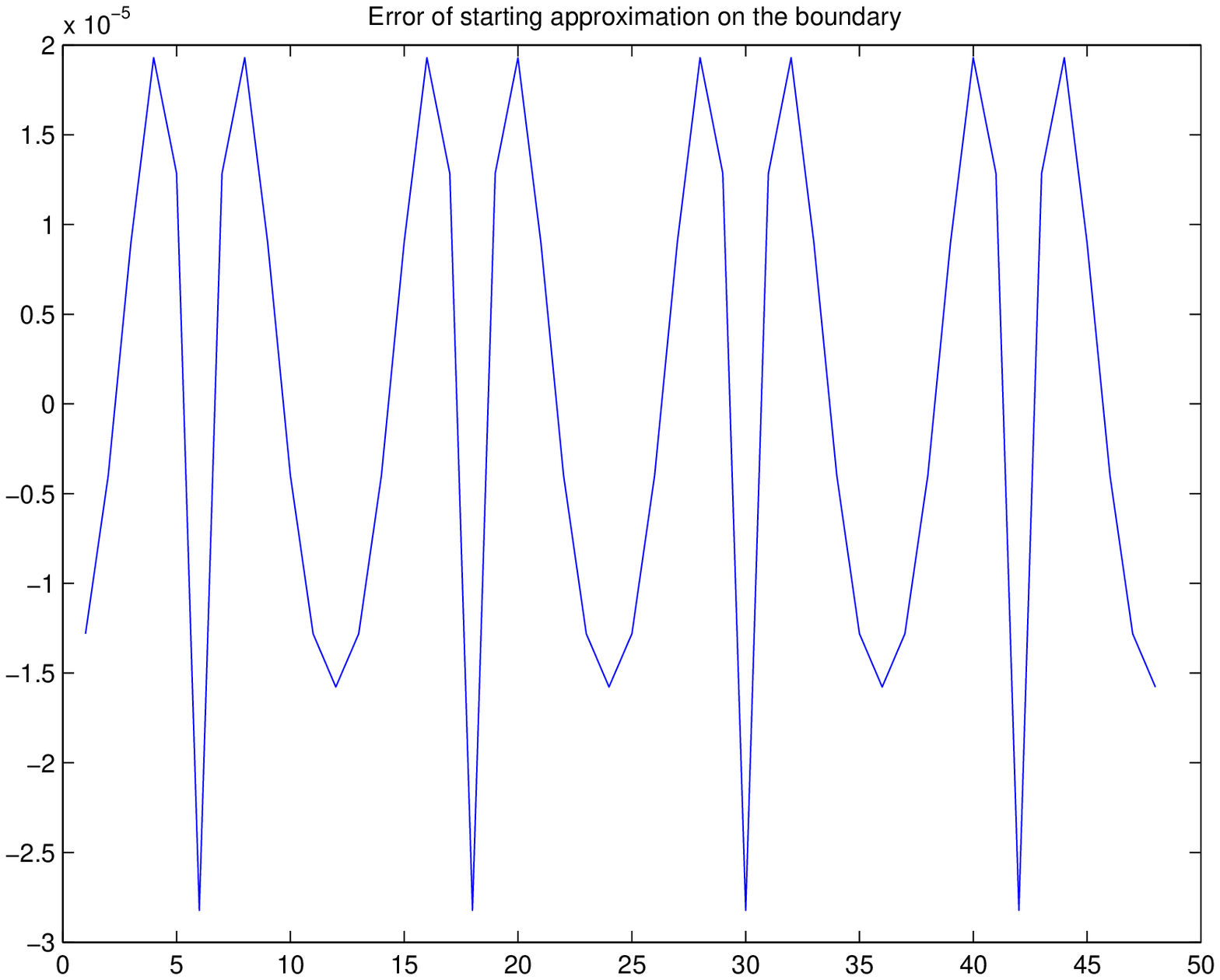}\\
\caption{Starting error on the boundary}\label{MA0112errbndstart}
\end{center}
\end{figure}
\begin{figure}[hbt]
\begin{center}
\includegraphics[width=10cm]{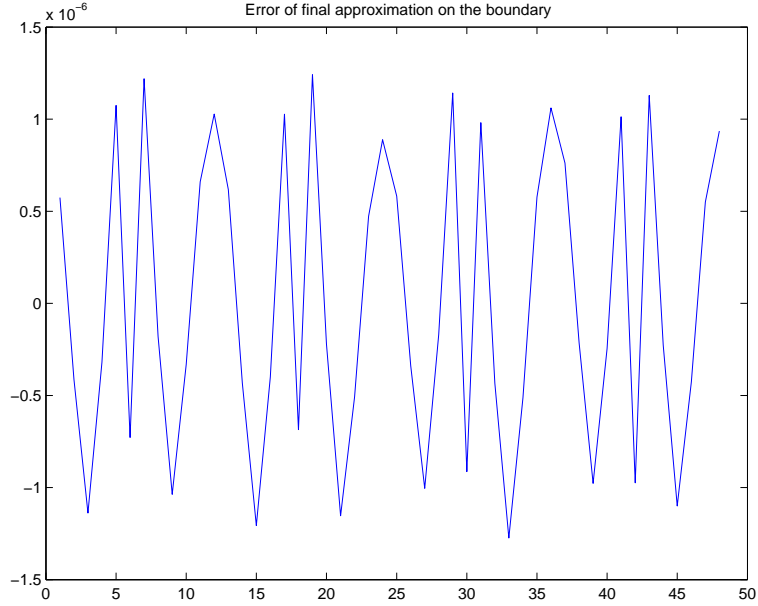}\\
\caption{Final error on the boundary}\label{MA0112errbnd}
\end{center}
\end{figure}
\begin{figure}[hbt]
\begin{center}
\includegraphics[width=10cm]{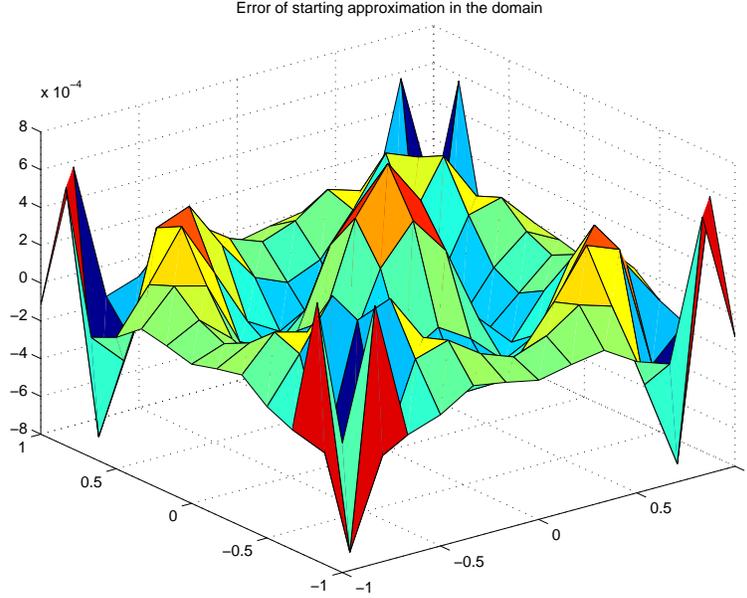}\\
\caption{Starting error in the PDE}\label{MA0112errintstart}
\end{center}
\end{figure}
\begin{figure}[hbt]
\begin{center}
\includegraphics[width=10cm]{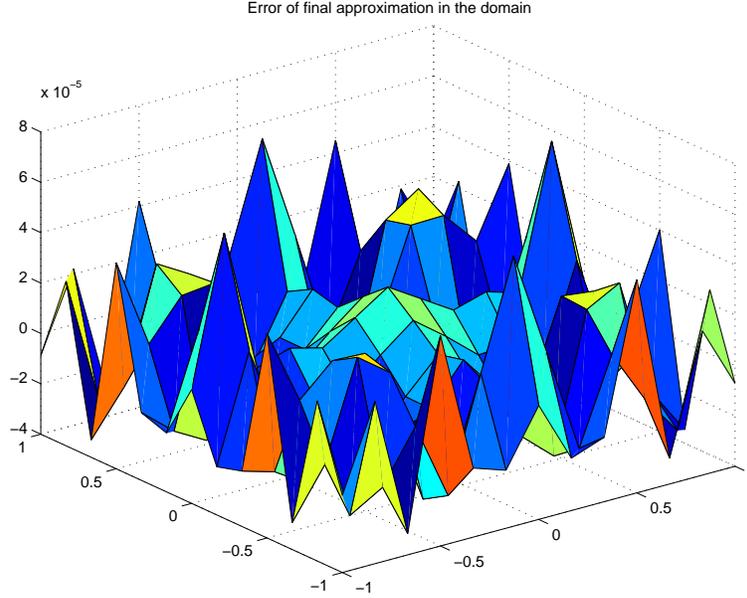}\\
\caption{Final error in the PDE}\label{MA0112errint}
\end{center}
\end{figure}
\begin{figure}[hbt]
\begin{center}
\includegraphics[width=10cm]{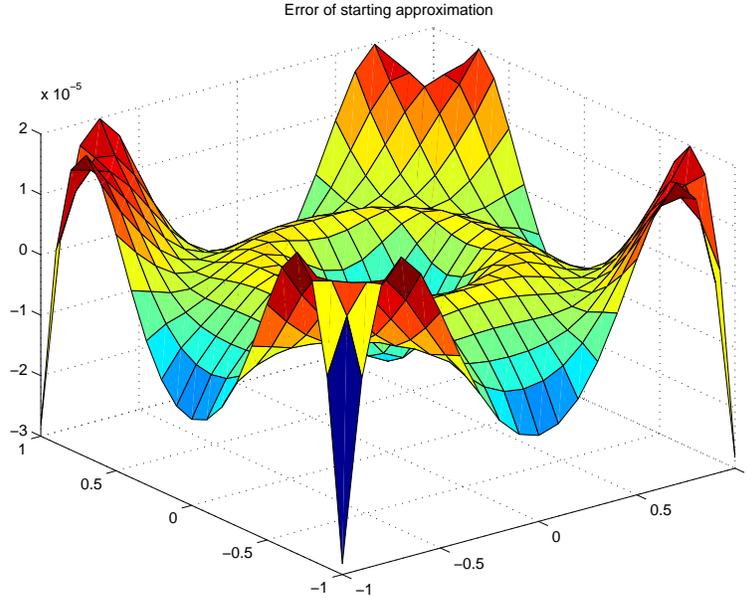}\\
\caption{Starting error wrt. the solution}\label{MA0112errstart}
\end{center}
\end{figure}
\begin{figure}[hbt]
\begin{center}
\includegraphics[width=10cm]{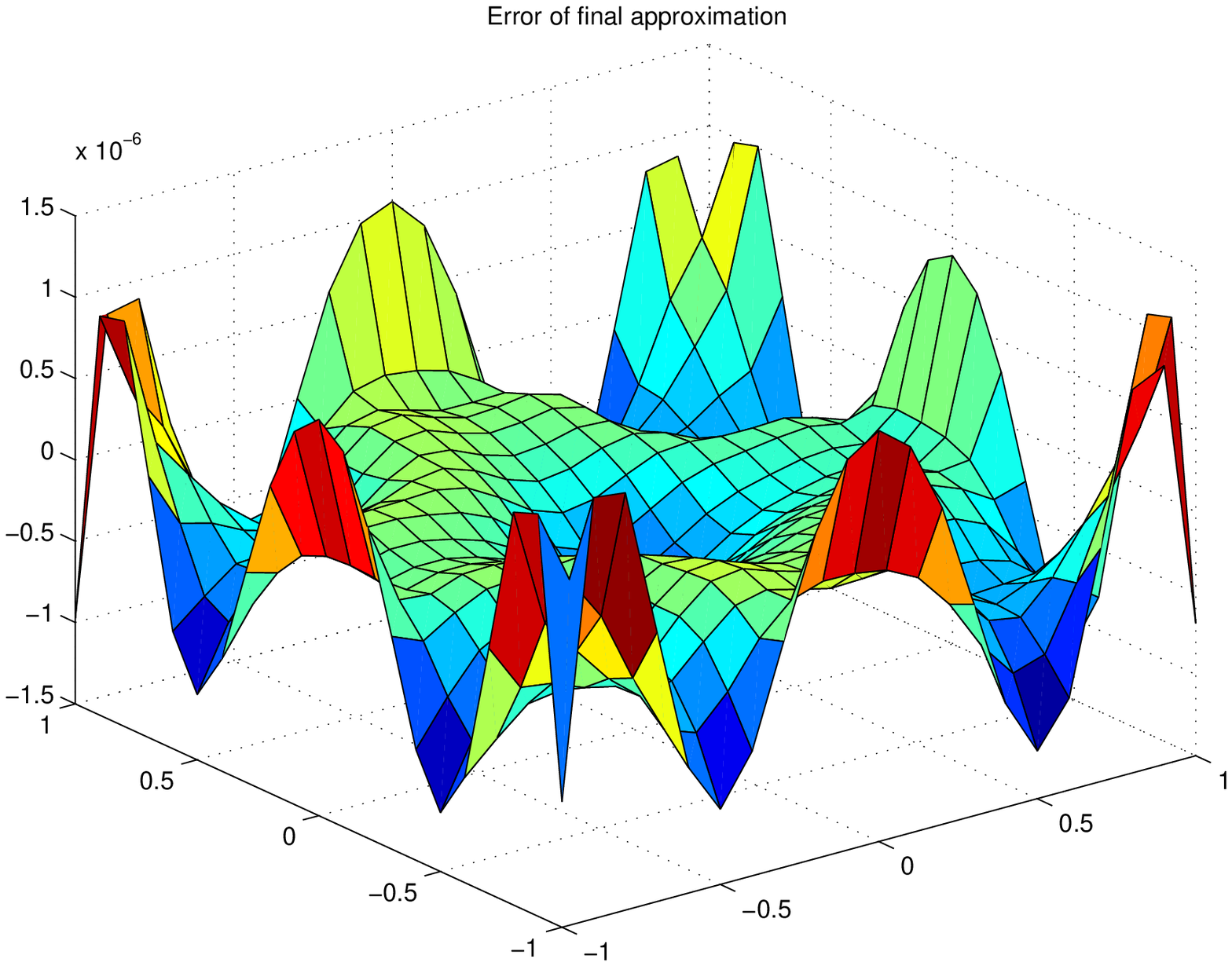}\\
\caption{Final error wrt. the solution}\label{MA0112err}
\end{center}
\end{figure}
\biglf
The results are rather promising and justify a more general
analysis to be provided in \cite{boehmer-schaback:2016-1}.
%****************************************************************
\section{MATLAB Implementation Details}\label{secMATDet}
%****************************************************************
% \subsection{Programming: Step 1}\label{SecP1}
The unknown coefficients are stored in a
triangular % $(M+1)\times(N+1)$
matrix $C$.
At each point $p_i=(x_i,y_i)$ we pre--calculate the
triangular Vandermonde matrix $Pu_i$ with entries $x_i^m\cdot y_i^n$ for 
$0\leq m+n\leq M$.
Since the coefficient matrix $C$ has the same shape,
we can take the element\-wise product
of $C$ and $Pu_i$, i.e. we
form $C.*Pu_i$ in MATLAB notation. Applying the {\tt sum} function of MATLAB
then yields $\MS(\MS(C.*Pu_i))=u^C(x_i,y_i)$. 
If $(x_i,y_i)$ is a boundary point, we form
$$
f_i(C):=f(x_i,y_i)-u^C(x_i,y_i)=f(x_i,y_i)-\MS(\MS(C.*Pu_i))
$$
and add $f_i(C)^2$ later into the minimization via
the MATLAB function {\tt nsqnonlin}. 
This cares for boundary value discretization.
\biglf
For PDE discretization, we work similarly. In particular,
the second derivatives
$$
\begin{array}{rcl}
u_{xx}(x,y)&=&\displaystyle{\sum_{m=2}^M\sum_{n=0}^{M-m}c_{mn}m(m-1)x^{m-2}y^n   },\\ 
u_{xy}(x,y)&=&\displaystyle{\sum_{m=1}^M\sum_{n=1}^{M-m}c_{mn}mnx^{m-1}y^{n-1}   },\\ 
u_{yy}(x,y)&=&\displaystyle{\sum_{m=0}^M\sum_{n=2}^{M-m}c_{mn}n(n-1)x^{m}y^{n-2}   }\\ 
\end{array} 
$$
are assembled at each point $(x_i,y_i)$ into triangular matrices
$ P_{xx}u_i,\,P_{xy}u_i,\,P_{yy}u_i $ with the above entries, 
%\begin{array}{rcl}
%P_{xx}u_i&:=&\left( m(m-1)x_i^{m-2}y_i^n  \right),\\
%P_{xy}u_i&:=&\left( mnx_i^{m-1}y_i^{n-1} \right),\\
%P_{yy}u_i&:=&\left( n(n-1)x_i^{m}y_i^{n-2} \right),\\
%\end{array}
%$$
and we get
$$
\begin{array}{rcl}
u_{xx}(x_i,y_i)(C) &=& \MS(\MS(C.* P_{xx}u_i)),\\
u_{xy}(x_i,y_i)(C) &=& \MS(\MS(C.* P_{xy}u_i)),\\
u_{yy}(x_i,y_i)(C) &=& \MS(\MS(C.* P_{yy}u_i))\\
\end{array}
$$
as matrix-valued functions of the coefficient matrix $C$.
If $(x_i,y_i)$ is a point where we
want the discretized PDE to be satisfied, we define
$$
\begin{array}{rcl}
f_i(C)&:=&u_{xx}^C(x_i,y_i)u_{yy}^C(x_i,y_i)-u_{xy}^C(x_i,y_i)^2-g(x_i,y_i)\\
&=&\MS(\MS(C.* P_{xx}u_i)).*
\MS(\MS(C.* P_{yy}u_i))\\
&& -(\MS(\MS(C.* P_{xy}u_i))).{^\wedge}2 -g(x_i,y_i)
\end{array} 
$$  
and add $f_i(C)^2$ into the minimization via {\tt nsqnonlin}.
This cares for the collocation of the nonlinear differential equation. 
\biglf
The total optimization problem then minimizes the sum of all these $f_i(C)^2$
with respect to the coefficient matrix $C$, using {\tt nlsqnonlin} of MATLAB. 
\biglf
An additional simplification turned triangular matrices into vectors.
Furthermore, the spatial discretization was changed with the polynomial
degree. For $M>0$ we used regular data with spacing  $h=2/M$, with $h=1$ for
$M=0$. Choosing finer spatial discretizations 
does not improve the results. To give boundary values more weight, we 
multiplied the boundary $f_i(C)$ with 10 throughout. 
%****************************************************************
 \section{Error Analysis}\label{SecConv}
% \red{RS test: \cite{boehmer-schaback:2016-1}, bitte woanders hinsetzen}\\
Since the paper \cite{boehmer-schaback:2013-1} reduces the convergence analysis
to the linearized strongly elliptic  problem under the above circumstances, we only need to deal
with the trial functions \eref{equSScxy} used for solving a standard elliptic
problem. Provided that sufficiently large sets of test points are chosen 
for each $M$,  
the paper \cite{schaback:2015-4} shows that the convergence rate 
is the one that arises when approximating the strong data of the solution,
i.e. the second derivatives of the solution by such polynomials. Thus the
convergence rate for the special case \eref{equexp} is exponential, as can be 
either taken from Bernstein-type theorems in Approximation Theory or directly
proven via
the exponentially convergent Taylor expansion. The numerical results
of Section \ref{SecR1} confirm this.
\biglf
It is much more difficult to assess how many test points are sufficient
for uniform stability.
Considering Theorem 5.1 and (5.3) of \cite{schaback:2015-4}, 
it suffices to prove conditions on the points that allow to bound $\|u^C\|_\infty $
uniformly in terms of the maximum or sum of the three discrete seminorms
$$
\begin{array}{rcl}
\|u^C\|_0&:=& \displaystyle{\max_{z_i^B\in Z^B}|u^C(z_i^B)|}\\
\|u^C\|_{xx}&:=& \displaystyle{ \max_{z_i^D\in Z^D}|u^C_{xx}(z_i^D)|}\\
\|u^C\|_{yy}&:=& \displaystyle{  \max_{z_i^D\in Z^D}|u^C_{yy}(z_i^D)|}\\
\end{array} 
$$
for arbitrary trial functions $u^C$ of the form \eref{equSScxy}.
Note that this is independent of PDEs. It is a problem of Approximation Theory.
We can ignore the mixed derivatives and the factors depending on $u^*$ here,
because uniform strong ellipticity allows to go over to the pure second
derivatives, at the expense of uniformly bounded factors depending on $u^*$.   
\biglf
Using the standard machinery for proofs of stability inequalities, as nicely summarized in
Chapter 3 of \cite{wendland:2005-1}, we start with fill distances $h_B$
and $h_D$ for the points on the boundary and the domain, respectively, and deal
with the boundary first. As a warm-up, consider the upper and lower boundary
lines, i.e. two lines of length 2
parallel to the $x$ axis. There, each $u^C$ is a univariate polynomial $p$
of degree at most $M$,
and we rewrite it as one on $[-1,+1]$. 
Now we know by Markov's inequality that 
$$
\|p'\|_{\infty,[-1,+1]}\leq M^2 \|p\|_{\infty,[-1,+1]}.
$$
For each point $x$ on the boundary lines there is a sample point $x_i$ at
distance at most $h_B$, and thus
$$
|p(x)|\leq |p(x_i)|+h_BM^2\|p\|_{\infty,[-1,+1]}
$$
proving $\|p\|_\infty\leq 2 \max_{x_i} |p(x_i)|$ on those boundary lines, if we
have $h_BM^2\leq 1/2$. This proves 
$$
\|u^C\|_{\infty,\Gamma}\leq 2 \max_{z_i^B\in Z_B}|u^C(z_i^B)|
$$ 
as soon as 
$$
h_B\leq \frac{1}{2}M^{-2}.
$$
Note that Approximation Theory
tells us that we can replace the exponent two by
one,
if ${\cal O}(M)$ points on the boundary are distributed in Chebyshev
style.
\biglf
The case with derivatives and interior points is not directly covered by the
standard theory. From \cite{braess:2001-1} we know that there is a
well--posedness inequality 
\bql{eqWPI}
\|u\|_{\infty,\overline{\Omega}}\leq \|u\|_{\infty,\Gamma}+C\|\Delta u\|_{\infty,\overline{\Omega}}
\eq
for uniformly elliptic problems, and this can be applied to trial functions. 
We already have the first term on the right--hand side under control and only
have to deal with the second term. Since all $\Delta u^C$ are polynomials as
well as $u^C$, we can use the standard logic as in  Chapter 3 of
\cite{wendland:2005-1} to get a bound of the form 
$$
\|\Delta u^C\|_{\infty,\Omega}\leq 2 \max_{z_i^D\in Z_D}|\Delta
u^C(z_i^D)|
\leq 2\|u^C\|_{xx}+ 2\|u^C\|_{yy}
$$
if $h_D\leq CM^{-2}$.
Now the theory in 
\cite{schaback:2010-2,boehmer-schaback:2013-1,schaback:2015-4}
implies
\begin{theorem}\label{TheColloc}
  The strong meshless
  collocation method for solving the Monge-Amp\`ere equation by trial
functions \eref{equSScxy} via sampling on test points with fill distances $h_B$
and $h_D$ on the boundary and the domain is uniformly stable if the fill
distances behave like ${\cal O}(M^{-2})$. Furthermore, the convergence for
$M\to\infty$ is exponential if the test points are chosen to guarantee
uniform stability via sufficient oversampling. \qed
\end{theorem} 
%%%%%%%%%%%%%%%%%%%%%%%%%%%%
%****************************************************************
\section{Summary and Outlook}\label{SecSO}
By a specific example, it was demonstrated
theoretically and numerically 
that a strong meshless discretization of the Monge-Amp\`ere equation
works successfully. 
This will generalize to subdomains of $[-1,+1]^2$ 
with Lipschitz boundaries 
on which bivariate polynomials are still polynomials.
The domains should have a
uniform interior cone condition.
\biglf
Likewise, other trial spaces can be used.
If the true solution is less smooth,
convergence rates will then be confined to how well second derivatives
of trial functions
approximate second derivatives of the true solution.
To guarantee
stability for sufficient oversampling, the trial spaces must
allow that the results of the previous section can be applied,  
\biglf
A much more thorough investigation
of meshless methods for fully nonlinear problems is in preparation
\cite{boehmer-schaback:2016-1}.
% \textbf{Klaus WHY TWICE ??? } RS: Weil es in der summary ist
% und manche nur Abstract und Summary lesen.
%%%%%%%%%%%%%%%%%%%%%%%%%%%%%%%%%%%%%%%%%%%%%%%%%%%%%%%%%
\bibliographystyle{plain}
%\bibliography{RSbib}

\begin{thebibliography}{10}

\bibitem{benamou-et-al:2010-1}
J.-D. Benamou, B.D. Froese, and A.M. Oberman.
\newblock Two numerical methods for the elliptic {M}onge-{A}mpère equation.
\newblock {\em ESAIM: Mathematical Modelling and Numerical Analysis},
  44(4):737--758, 2010.

\bibitem{boehmer:2004-1}
K.~B{\"{o}}hmer.
\newblock On finite element methods for fully nonlinear elliptic equations of
  second order.
\newblock {\em SIAM J. Numer. Anal.}, 3:1212--1249, 2008.

\bibitem{boehmer:2010-1}
K.~B\"ohmer.
\newblock {\em Numerical Methods for Nonlinear Elliptic Differential Equations,
  a Synopsis}.
\newblock Oxford University Press, Oxford, 2010.

\bibitem{boehmer-schaback:2013-1}
K.~B\"ohmer and R.~Schaback.
\newblock A nonlinear discretization theory.
\newblock {\em Journal of Computational and Applied Mathematics}, 254:204--219,
  2013.

\bibitem{boehmer-schaback:2016-1}
K.~B\"ohmer and R.~Schaback.
\newblock Nonlinear discretization theory applied to meshfree methods and the
  {M}onge-{A}mp\`{e}re equation.
\newblock Fachbereich {M}athematik und {I}nformatik,
  {P}hilipps--{U}niversit{\"{a}}t {M}arburg, in preparation, 2016.

\bibitem{braess:2001-1}
D.~Braess.
\newblock {\em Finite Elements. Theory, Fast Solvers and Applications in Solid
  Mechanics}.
\newblock Cambridge University Press, 2001.
\newblock Second edition.

\bibitem{brenner-et-al:2011-1}
S.C. Brenner, T.~Gudi, M.~Neilan, and L.Y. Sung.
\newblock {C}$^0$ penalty methods for the fully nonlinear {M}onge-{A}mp\`ere
  equation.
\newblock {\em Mathematics of Computation}, 80:1979 -- 1995, 2011.

\bibitem{davydov:2006-1}
O.~Davydov.
\newblock Smooth finite elements and stable splitting.
\newblock Fachbereich {M}athematik und {I}nformatik,
  {P}hilipps--{U}niversit{\"{a}}t {M}arburg, 2007.

\bibitem{davydov-saaed:2013-1}
O.~Davydov and A.~Saeed.
\newblock Numerical solution of fully nonlinear elliptic equations by
  {B}{\"o}hmer's method.
\newblock {\em J. Comp.Appl.Math.}, 254:43--54, 2013.

\bibitem{fasshauer-mccourt:2015-1}
G.~Fasshauer and M.~McCourt.
\newblock {\em Kernel-based Approximation Methods using MATLAB}, volume~19 of
  {\em Interdisciplinary Mathematical Sciences}.
\newblock World Scientific, Singapore, 2015.

\bibitem{feng-neilan:2009-1}
X.~Feng and M.~Neilan.
\newblock Mixed finite element methods for the fully nonlinear
  {M}onge-{A}mp\`ere equation based on the vanishing moment method.
\newblock {\em SIAM Journal on Numerical Analysis}, 47:1226 -- 1250, 2009.

\bibitem{froese-oberman:2011-1}
B.D. Froese and A.M. Oberman.
\newblock Convergent finite difference solvers for viscosity solutions of the
  elliptic {M}onge-{A}mp\`ere equation in dimensions two and higher.
\newblock {\em SIAM Journal on Numerical Analysis}, 49(4):1692--1714, 2011.

\bibitem{froese-oberman:2013-1}
B.D. Froese and A.M. Oberman.
\newblock Convergent filtered schemes for the {M}onge--{A}mp\`ere partial
  differential equation.
\newblock {\em SIAM Journal on Numerical Analysis}, 51(1):423--444, 2013.

\bibitem{li-liu:2017-1}
Q.~Li and Z.Y. Liu.
\newblock Solving the 2-{D} elliptic {M}onge-{A}mp{\`e}re equation by a
  {K}ansa's method.
\newblock {\em Acta Mathematicae Applicatae Sinica, English Series},
  33(2):269--276, Apr 2017.

\bibitem{liu-et-al:2017-1}
J.~Liu, B.D. Froese, A.M. Oberman, and M.Q. Xiao.
\newblock A multigrid scheme for 3{D} {M}onge-{A}mpère equations.
\newblock {\em International Journal of Computer Mathematics},
  94(9):1850--1866, 2017.

\bibitem{liu-he:2013-1}
Z.Y. Liu and Y.~He.
\newblock Cascadic meshfree method for the elliptic {M}onge-{A}mpère equation.
\newblock {\em Engineering Analysis with Boundary Elements}, 37(7):990 -- 996,
  2013.

\bibitem{liu-he:2014-1}
Z.Y. Liu and Y.~He.
\newblock An iterative meshfree method for the elliptic {M}onge-{A}mp\`ere
  equation in 2{D}.
\newblock {\em Numerical Methods for Partial Differential Equations},
  30(5):1507--1517, 2014.

\bibitem{oberman:2008-1}
A.~Oberman.
\newblock Wide stencil finite difference schemes for the elliptic
  {M}onge-{A}mp\`ere equations and functions of the eigenvalues of the
  {H}essian.
\newblock {\em Discrete Contin. Dyn. Syst. Ser B 10(1)}, pages 221--238, 2008.

\bibitem{schaback:2010-2}
R.~Schaback.
\newblock Unsymmetric meshless methods for operator equations.
\newblock {\em Numerische Mathematik}, 114:629--651, 2010.

\bibitem{schaback:2015-4}
R.~Schaback.
\newblock All well--posed problems have uniformly stable and convergent
  discretizations.
\newblock {\em Numerische Mathematik}, 132:597--630, 2015.

\bibitem{wendland:2005-1}
H.~Wendland.
\newblock {\em Scattered Data Approximation}.
\newblock Cambridge University Press, 2005.

\end{thebibliography}

\end{document}